\documentclass[12pt]{article}
\usepackage{latexsym,amssymb,amsmath,enumerate,geometry,float,graphicx}
\usepackage{graphics}

\newtheorem{definition}{Definition}
\newtheorem{theorem}[definition]{Theorem}

\newtheorem{lemma}[definition]{Lemma}

\usepackage{lineno}

\begin{document}

\title{Independence in Uniform Linear Triangle-free Hypergraphs}
\author{Piotr Borowiecki$^1$, Michael Gentner$^2$, Christian L\"{o}wenstein$^2$, Dieter Rautenbach$^2$}

%\linenumbers

\date{}

\maketitle

\begin{center}
{\small 
$^1$ 
Department of Algorithms and System Modeling, 
Faculty of Electronics, Telecommunications and Informatics,
Gda\'{n}sk University of Technology,
80-233 Gda\'{n}sk, Poland,
email: pborowie@eti.pg.gda.pl\\[3mm]
$^2$
Institute of Optimization and Operations Research, 
Ulm University, 
D-89069 Ulm,
Germany,\\
email: $\{$michael.gentner, christian.loewenstein, dieter.rautenbach$\}$@uni-ulm.de
}
\end{center}

\begin{abstract}
The independence number $\alpha(H)$ of a hypergraph $H$ 
is the maximum cardinality of a set of vertices of $H$ that does not contain an edge of $H$.
Generalizing Shearer's classical lower bound on the independence number of triangle-free graphs
(J. Comb. Theory, Ser. B 53 (1991) 300-307),
and considerably improving recent results of
Li and Zang (SIAM J. Discrete Math. 20 (2006) 96-104)
and 
Chishti et al.~(Acta Univ. Sapientiae, Informatica 6 (2014) 132-158), 
we show that 
$$\alpha(H)\geq \sum_{u\in V(H)}f_r(d_H(u))$$
for an $r$-uniform linear triangle-free hypergraph $H$ with $r\geq 2$, where 
\begin{eqnarray*}
f_r(0)&=&1\mbox{, and }\\
f_r(d)&=&\frac{1+\Big((r-1)d^2-d\Big)f_r(d-1)}{1+(r-1)d^2}\mbox{ for $d\geq 1$.}
\end{eqnarray*}
\end{abstract}



%\pagebreak

\section{Introduction}

We consider finite {\it hypergraphs} $H$, which are ordered pairs $(V(H),E(H))$ of two sets, where 
$V(H)$ is the finite set of {\it vertices of $H$} and
$E(H)$ is the {\it set of edges of $H$}, which are subsets of $V(H)$.
The {\it order $n(H)$ of $H$} is the cardinality of $V(H)$.
The {\it degree $d_H(u)$ of a vertex $u$ of $H$} is the number of edges of $H$ that contain $u$.
The {\it average degree $d(H)$ of $H$} is the arithmetic mean of the degrees of its vertices.
Two distinct vertices of $H$ are {\it adjacent} or {\it neighbors} if some edge of $H$ contains both.
The {\it neighborhood $N_H(u)$ of a vertex $u$ of $H$} is the set of vertices of $H$ that are adjacent to $u$.
For a set $X$ of vertices of $H$, the hypergraph $H-X$
arises from $H$ by removing from $V(H)$ all vertices in $X$
and removing from $E(H)$ all edges that intersect $X$.
If every two distinct edges of $H$ share at most one vertex, then $H$ is {\it linear}.
If $H$ is linear and for every two distinct non-adjacent vertices $u$ and $v$ of $H$, 
every edge of $H$ that contains $u$ contains at most one neighbor of $v$, then $H$ is {\it double linear}.
If there are not three distinct vertices $u_1$, $u_2$, and $u_3$ of $H$ and three distinct edges $e_1$, $e_2$, and $e_3$ of $H$ 
such that $\{ u_1,u_2,u_3\}\setminus \{ u_i\}\subseteq e_i$ for $i\in \{ 1,2,3\}$, then $H$ is {\it triangle-free}.
A set $I$ of vertices of $H$ is a {\it (weak) independent set of $H$} if no edge of $H$ is contained in $I$.
The {\it (weak) independence number $\alpha(H)$ of $H$} is the maximum cardinality of an independent set of $H$.
If all edges of $H$ have cardinality $r$, then $H$ is {\it $r$-uniform}.
If $H$ is $2$-uniform, then $H$ is referred to as a {\it graph}.

The independence number of (hyper)graphs is a well studied computationally hard parameter.
Caro \cite{ca} and Wei \cite{we} proved a classical lower bound on the independence number of graphs,
which was extended to hypergraphs by Caro and Tuza \cite{catu}.
Specifically, for an $r$-uniform hypergraph $H$, Caro and Tuza \cite{catu} proved
$$\alpha(H)\geq \sum_{u\in V(H)}f_{CT(r)}(d_H(u)),$$ 
where 
$f_{CT(r)}(d)={d+\frac{1}{r-1}\choose d}^{-1}$.
Thiele \cite{th} generalized Caro and Tuza's bound to general hypergraphs;
see \cite{bogohara} for a very simple probabilistic proof of Thiele's bound.
Originally motivated by Ramsey theory, 
Ajtai et al.~\cite{ajkosz} showed that $\alpha(G)=\Omega\left(\frac{\ln d(G)}{d(G)}n(G)\right)$ for every triangle-free graph $G$.
Confirming a conjecture from \cite{ajkosz} concerning the implicit constant, 
Shearer \cite{sh1} improved this bound to $\alpha(H)\geq f_{S_1}(d(G))n(G)$, 
where $f_{S_1}(d)=\frac{d\ln d-d+1}{(d-1)^2}$.
In \cite{sh1} the function $f_{S_1}$ arises as a solution of the differential equation
$$(d+1)f(d) = 1+(d-d^2)f'(d)\mbox{ and }f(0)=1.$$
In \cite{sh2} Shearer showed that 
$$\alpha(G)\geq \sum_{u\in V(G)}f_{S_2}(d_G(u))$$
for every triangle-free graph $G$,
where $f_{S_2}$ solves the difference equation
$$(d+1)f(d) = 1+(d-d^2)\Big(f(d)-f(d-1)\Big)\mbox{ and }f(0)=1.$$
Since $f_{S_1}(d)\leq f_{S_2}(d)$ for every non-negative integer $d$, and $f_{S_1}$ is convex,
Shearer's bound from \cite{sh2} is stronger than his bound from \cite{sh1}.

Li and Zang \cite{liza} adapted Shearer's approach to hypergraphs and obtained the following.

\begin{theorem}[Li and Zang \cite{liza}]\label{theoremlizang}
Let $r$ and $m$ be positive integers with $r\geq 2$.

If $H$ is an $r$-uniform double linear hypergraph 
such that the maximum degree of every subhypergraph of $H$ induced by the neighborhood of a vertex of $H$ is less than $m$,
then 
$$\alpha(H)\geq \sum_{u\in V(H)}f_{LZ(r,m)}(d_H(u)),$$ 
where 
\begin{eqnarray*}
f_{LZ(r,m)}(x)&=&\frac{m}{B}\int_0^1\frac{(1-t)^{\frac{a}{m}}}{t^b(m-(x-m)t)}dt,
\end{eqnarray*}
$a=\frac{1}{(r-1)^2}$,
$b=\frac{r-2}{r-1}$, and
$B=\int_0^1(1-t)^{\left(\frac{a}{m}-1\right)}t^{-b}dt$.
\end{theorem}
Note that for $r\geq 2$, 
an $r$-uniform linear hypergraph $H$ is triangle-free if and only if 
it is double linear and the maximum degree of every subhypergraph of $H$ induced by the neighborhood of a vertex of $H$ is less than $1$.
Therefore, since $f_{S_1}=f_{LZ(2,1)}$ and $f_{S_1}$ is convex, 
Theorem \ref{theoremlizang} implies Shearer's bound from \cite{sh1}.
Nevertheless, 
since $f_{S_1}(d)<f_{S_2}(d)$ for every integer $d$ with $d\geq 2$,
Shearer's bound from \cite{sh2} does not quite follow from Theorem \ref{theoremlizang}.

In \cite{chzhpiiv} Chishti et al.~presented another version of Shearer's bound from \cite{sh1} for hypergraphs.

\begin{theorem}[Chishti et al.~\cite{chzhpiiv}]\label{theoremchzhpiiv}
Let $r$ be an integer with $r\geq 2$.

If $H$ is an $r$-uniform linear triangle-free hypergraph, then 
$$\alpha(H)\geq f_{CZPI(r)}(d(H)) n(H),$$ 
where 
\begin{eqnarray*}
f_{CZPI(r)}(x)&=&\frac{1}{r-1}\int_0^1\frac{1-t}{t^b(1-((r-1)x-1)t)}dt
\end{eqnarray*}
and $b=\frac{r-2}{r-1}$.
\end{theorem}
Since $f_{S_1}=f_{CZPI(2)}$, for $r=2$, the last result coincides with Shearer's bound from \cite{sh1}.

\begin{figure}[H]
\begin{center}
$\mbox{}$
\hfill
\includegraphics[width=0.45\textwidth]{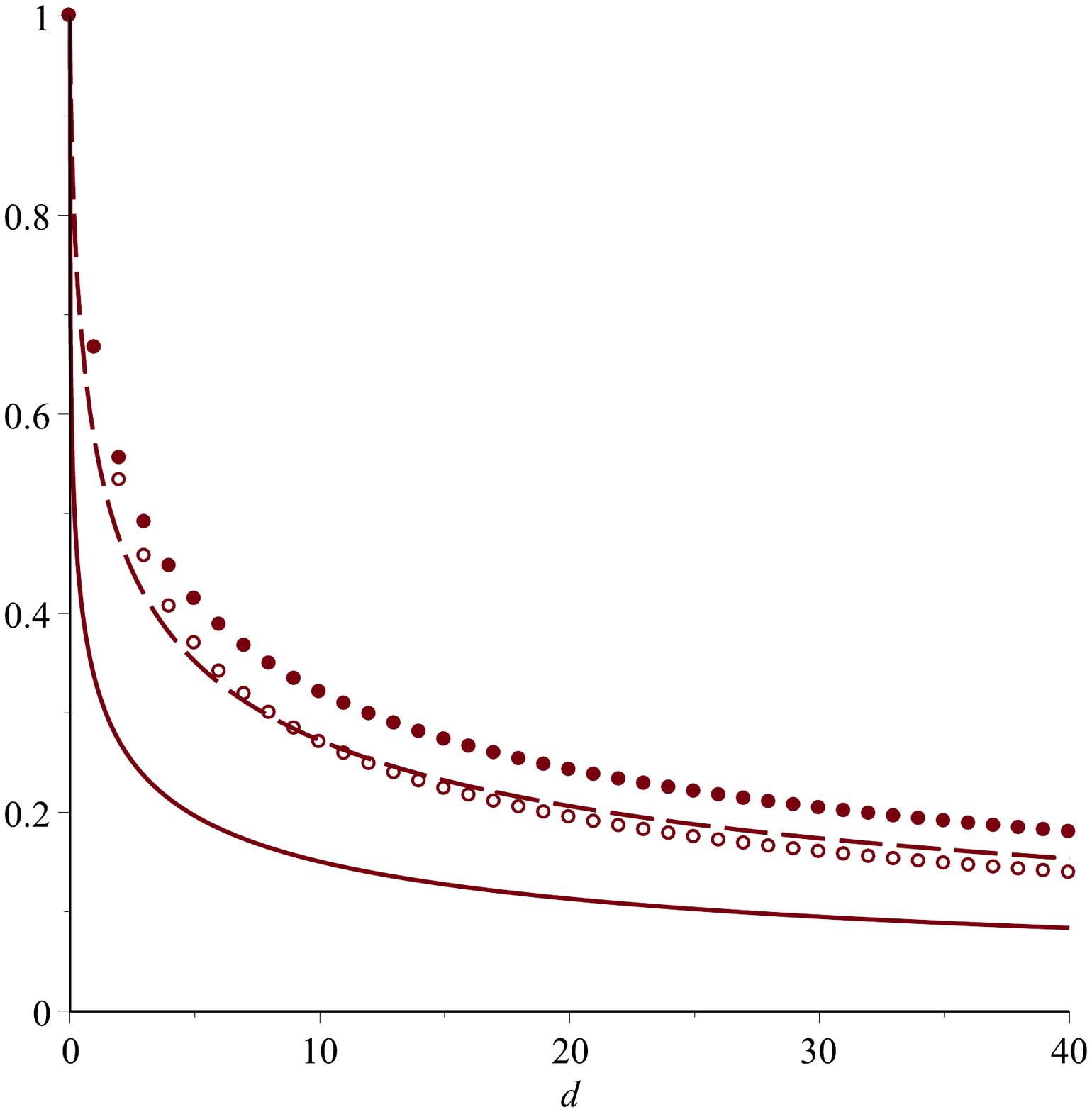}
\hfill
\includegraphics[width=0.45\textwidth]{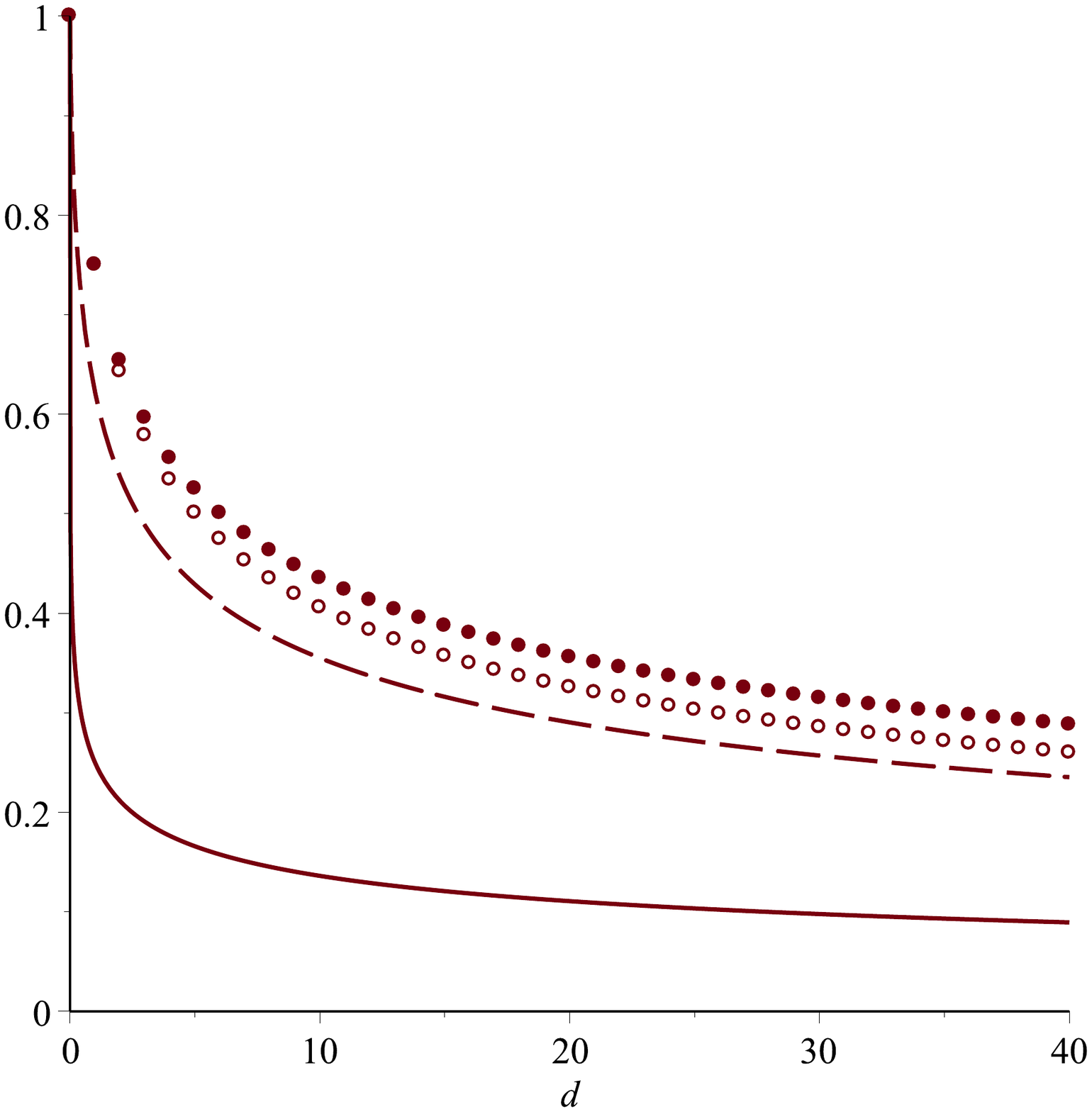}
\hfill
$\mbox{}$
\end{center}
\vspace{-2cm}
\caption{The values of  
$f_{LZ(r,1)}(d)$ (line), 
$f_{CZPI(r)}(d)$ (dashed line), 
$f_{CT(r)}(d)$ (empty circles), and
$f_r(d)$ (solid circles)
for $0\leq d\leq 40$ and
$r=3$ (left) and
$r=4$ (right).}\label{fig0}
\end{figure}

A drawback of the bounds in Theorem \ref{theoremlizang} and Theorem \ref{theoremchzhpiiv} is 
that they are very often weaker than Caro and Tuza's bound \cite{catu}, 
which holds for a more general class of hypergraphs.
See Figure \ref{fig0} for an illustration.

\vspace{-1.5cm}

\noindent In the present paper we extend Shearer's approach from \cite{sh2} and 
establish a lower bound on the independence number of a uniform linear triangle-free hypergraph
that considerably improves Theorem \ref{theoremlizang} and Theorem \ref{theoremchzhpiiv}
and is systematically better than Caro and Tuza's bound.

For further related results we refer to 
Ajtai et al. \cite{AKPSS},
Duke et al. \cite{DLR},
Dutta et al.~\cite{dumusu} and
Kostochka et al.~\cite{komuve}.
Note that our main result provides explicit values when applied to a specific hypergraph
but that we do not completely understand its asymptotics.
In contrast to that, results as in \cite{AKPSS,DLR,dumusu} 
are essentially asymptotic statements 
but are of limited value when applied to a specific hypergraph.

\section{Results}

For an integer $r$ with $r\geq 2$, let $f_r:\mathbb{N}_0\to\mathbb{R}_0$ be such that 
\begin{eqnarray*}
f_r(0)&=&1\mbox{ and}\\
f_r(d)&=&\frac{1+\Big((r-1)d^2-d\Big)f_r(d-1)}{1+(r-1)d^2}
\end{eqnarray*}
for every positive integer $d$.

\begin{lemma}\label{lemma2}
If $r$ and $d$ are integers with $r\geq 2$ and $d\geq 0$, then
$f_r(d)-f_r(d+1)\geq f_r(d+1)-f_r(d+2)$.
\end{lemma}
{\it Proof:} Substituting within the inequality
$f_r(d)-2f_r(d+1)+f_r(d+2)\geq 0$
first 
$f_r(d+2)$ with 
$$\frac{1+\Big((r-1)(d+2)^2-(d+2)\Big)f_r(d+1)}{1+(r-1)(d+2)^2}$$
and then 
$f_r(d+1)$ with 
$$\frac{1+\Big((r-1)(d+1)^2-(d+1)\Big)f_r(d)}{1+(r-1)(d+1)^2},$$
and solving it for $f_r(d)$, 
it is straightforward but tedious to verify that it is equivalent to $f_r(d)\geq L(r,d)$
where 
$$L(r,d)=\frac{(2r-1)d+3r}{r(d^2+5d+5)}.$$
Therefore, in order to complete the proof, it suffices to show $f_r(d)\geq L(r,d)$.
For $d=0$, we have $f_r(0)=1>\frac{3}{5}=L(r,0)$.
Now, let $f(d)\geq L(r,d)$ for some non-negative integer $d$.
Since $(r-1)(d+1)^2-(d+1)\geq 0$,
we obtain by a straightforward yet tedious calculation
\begin{eqnarray*}
f(d+1)-L(r,d+1)
&=&\frac{1+\Big((r-1)(d+1)^2-(d+1)\Big)f(d)}{1+(r-1)(d+1)^2}-L(r,d+1)\\
&\geq&\frac{(1+\Big((r-1)(d+1)^2-(d+1)\Big)L(r,d)}{1+(r-1)(d+1)^2}-L(r,d+1)\\
&=& \frac{2(1+(r-1)(d+2)^2)}{r(d^2+7d+11)(d^2+5d+5)},
\end{eqnarray*}
which is positive for $r\geq 2$.
Therefore, $f(d+1)\geq L(r,d+1)$,
which completes the proof by an inductive argument.
$\Box$

\bigskip

\noindent The following is our main result.

\begin{theorem}\label{theorem1}
Let $r$ be an integer with $r\geq 2$.

If $H$ is an $r$-uniform linear triangle-free hypergraph, then 
$$\alpha(H)\geq \sum_{u\in V(H)}f_r(d_H(u)).$$
\end{theorem}

\bigskip

\noindent Before we proceed to the proof,
we compare our bound 
to the bounds
of Caro and Tuza \cite{catu},
Li and Zang \cite{liza}, and
Chishti et al.~\cite{chzhpiiv}.
Figure \ref{fig0} illustrates some specific values. 
An inspection of Li and Zang's proof in \cite{liza} reveals that 
they actually prove a lower bound on the so-called {\it strong independence number},
which is defined as the maximum cardinality of a set of vertices that does not contain two adjacent vertices. 
Therefore, especially for large values of $r$, 
Theorem \ref{theoremlizang} is much weaker than Theorem \ref{theoremchzhpiiv}.
In fact, it is quite natural that it is worse by a factor of about $r-1$.

As we show now,
our bound is systematically better than Caro and Tuza's bound \cite{catu}.

\pagebreak

\begin{lemma}\label{lemma3}
If $r$ and $d$ are integers with $r\geq 3$ and $d\geq 2$, 
then $f_r(d)>f_{CT(r)}(d)$.
\end{lemma}
{\it Proof:} 
Note that  $f_r(0)=f_{CT(r)}(0)=1$, $f_r(1)=f_{CT(r)}(1)=\frac{r-1}{r}$, and $f_{CT(r)}(d)=\frac{d}{d+\frac{1}{r-1}}f_{CT(r)}(d-1)$ for $d\in \mathbb{N}$,
which immediately implies that $f_{CT(r)}(d)<\frac{r-1}{r}$ for $d\geq 2$.
Now, if $f_r(d-1)\geq f_{CT(r)}(d-1)$ for some $d\geq 2$, then
\begin{eqnarray*}
f_r(d)-f_{CT(r)}(d) 
&=& 
\frac{1+\Big((r-1)d^2-d\Big)f_r(d-1)}{1+(r-1)d^2}-f_{CT(r)}(d)\\
&\geq & 
\frac{1+\Big((r-1)d^2-d\Big)f_{CT(r)}(d-1)}{1+(r-1)d^2}-f_{CT(r)}(d)\\
&=&
\frac{1+\Big((r-1)d^2-d\Big)\frac{1+(r-1)d}{(r-1)d}f_{CT(r)}(d)}{1+(r-1)d^2}-f_{CT(r)}(d)\\
&=&
\frac{1-\frac{r}{r-1}f_{CT(r)}(d)}{1+(r-1)d^2}\\
&>& 0,
\end{eqnarray*}
that is, $f_r(d)>f_{CT(r)}(d)$,
which completes the proof by an inductive argument. $\Box$

\bigskip

\noindent For $r=2$, Lemma \ref{lemma3} would state that Shearer's bound \cite{sh2}
is better than Caro \cite{ca} and Wei's bound \cite{we}, 
which is known.

We proceed to the proof of Theorem \ref{theorem1}.

\bigskip

\noindent {\it Proof of Theorem \ref{theorem1}:}
We prove the statement by induction on $n(H)$.
If $H$ has no edge, then $\alpha(H)=n(H)$, which implies the desired result for $n(H)\leq r-1$.
Now let $n(H)\geq r$.
If $H$ has a vertex $x$ with $d_H(x)=0$, 
then $f_r(d_H(x))=1$ and, by induction,
$$\alpha(H)
\geq 1+\alpha(H-x)\geq f_r(d_H(x))+\sum_{u\in V(H)\setminus \{ x\}}f_r(d_{H-x}(u))
=\sum_{u\in V(H)}f_r(d_H(u)).$$
Hence we may assume that $H$ has no vertex of degree $0$.

Since $H$ is $r$-uniform and linear, for every two edges $e_1$ and $e_2$ with $e_1\cap e_2=\{ u\}$ for some vertex $u$ of $H$,
the sets $e_1\setminus \{ u\}$ and $e_2\setminus \{ u\}$ are disjoint and of order $r-1$.
Therefore, for every vertex $u$ of $H$, there is a set ${\cal R}(u)$ of $r-1$ sets of neighbors of $u$ such that 
every neighbor of $u$ belongs to exactly one of the sets in ${\cal R}(u)$, and
$|e\cap R|=1$ for every edge $e$ of $H$ with $u\in e$ and every $R\in {\cal R}(u)$.

If $x$ is a vertex of $H$ and $R\in {\cal R}(x)$ is such that 
$$1+\sum_{u\in V(H)\setminus (\{ x\}\cup R)}f_r(d_{H-(\{ x\}\cup R)}(u))\geq \sum_{u\in V(H)}f_r(d_H(u)),$$
then the statement follows by induction, 
because $\alpha(H)\geq 1+\alpha(H-(\{ x\}\cup R))$.
Therefore, in order to complete the proof, it suffices to show that 
the following term is non-negative:
\begin{eqnarray*}
P&=&\sum_{x\in V(H)}\sum_{R\in {\cal R}(x)}
\left(1+\sum_{u\in V(H)\setminus (\{ x\}\cup R)}f_r(d_{H-(\{ x\}\cup R)}(u))-\sum_{u\in V(H)}f_r(d_H(u))\right).
\end{eqnarray*}
Since $H$ is linear and triangle-free, we have $d_{H-(\{ x\}\cup R)}(z)=d_H(z)-|N_H(z)\cap R|$
for every vertex $z$ in $V(H)\setminus (\{ x\}\cup R)$.
Trivially, $d_{H-(\{ x\}\cup R)}(z)=d_H(z)$ for $z\not\in N_H(R)$, and hence $P$ equals $P_1+P_2$, where
\begin{eqnarray*}
P_1 &=& \sum_{x\in V(H)}\sum_{R\in {\cal R}(x)}\left(1-f_r(d_H(x))-\sum_{y\in R}f_r(d_H(y))\right)\mbox{ and }\\
P_2 &=& \sum_{x\in V(H)}\sum_{R\in {\cal R}(x)}\sum_{z\in N_H(R)\setminus \{ x\}}\Big(f_r(d_H(z)-|N_H(z)\cap R|)-f_r(d_H(z))\Big)
\end{eqnarray*}
Since for every vertex $u$ of $H$, there are exactly $(r-1)d_H(u)$ many vertices $v$ of $H$
such that $u$ belongs to exactly one of the sets in ${\cal R}(v)$, we have 
\begin{eqnarray*}
P_1 &= & \sum_{x\in V(H)}\Big((r-1)-(r-1)(d_H(x)+1)f_r(d_H(x))\Big).
\end{eqnarray*}
Since $f_r(d-1)-f_r(d)$ is decreasing by Lemma \ref{lemma2}, we have
$f_r(d-n)-f_r(d)\geq n(f_r(d-1)-f_r(d))$ for all positive integers $d$ and $n$ with $n<d$.
Therefore, 
\begin{eqnarray*}
P_2 &\geq&
\sum_{x\in V(H)}\sum_{R\in {\cal R}(x)}\sum_{z\in N_H(R)\setminus \{ x\}}
|N_H(z)\cap R|\Big(f_r(d_H(z)-1)-f_r(d_H(z))\Big)\\
&=& 
\sum_{x\in V(H)}\sum_{R\in {\cal R}(x)}\sum_{z\in N_H(R)\setminus \{ x\}}
\sum_{y\in R}|N_H(z)\cap \{ y\}|\Big(f_r(d_H(z)-1)-f_r(d_H(z))\Big)\\
&=& 
\sum_{x\in V(H)}\sum_{R\in {\cal R}(x)}\sum_{y\in R}\sum_{z\in N_H(R)\setminus \{ x\}}
|N_H(z)\cap \{ y\}|\Big(f_r(d_H(z)-1)-f_r(d_H(z))\Big)\\
&=& 
\sum_{x\in V(H)}\sum_{R\in {\cal R}(x)}\sum_{y\in R}\sum_{z\in N_H(y)\setminus \{ x\}}
\Big(f_r(d_H(z)-1)-f_r(d_H(z))\Big).
\end{eqnarray*}
Let $T$ be the set of all $4$-tupels $(x,R,y,z)$ with 
$x\in V(H)$,
$R\in {\cal R}(x)$,
$y\in R$, and 
$z\in N_H(y)\setminus \{ x\}$.
Note that $y\in N_H(z)$ for every $(x,R,y,z)$ in $T$.
Since $H$ is linear, 
for a given vertex $z$ of $H$ and a given neighbor $y$ of $z$,
there are $(r-1)d_H(y)-1$ many vertices $x$ of $H$ 
with $y\in R$ for some $R$ in ${\cal R}(x)$ 
and $z\in N_H(y)\setminus \{ x\}$.
Furthermore, by the properties of ${\cal R}(x)$, 
given $x$ and $y$, the set $R$ in ${\cal R}(x)$ with $y\in R$ is unique.
Therefore, 
\begin{eqnarray*}
P_2&\geq &
\sum_{x\in V(H)}\sum_{R\in {\cal R}(x)}\sum_{y\in R}\sum_{z\in N_H(y)\setminus \{ x\}}\Big(f_r(d_H(z)-1)-f_r(d_H(z))\Big)\\
&=&
\sum_{z\in V(H)}
\sum_{y\in N_H(z)}\Big((r-1)d_H(y)-1\Big)\Big(f_r(d_H(z)-1)-f_r(d_H(z))\Big).
\end{eqnarray*}
Let ${\cal E}$ be the edge set of the graph that arises from $H$ by replacing every edge of $H$ by a clique,
that is, ${\cal E}$ is the set of all sets containing exactly two adjacent vertices of $H$.

We obtain 
\begin{eqnarray*}
P_2&\geq &
\sum_{z\in V(H)}
\sum_{y\in N_H(z)}\Big((r-1)d_H(y)-1\Big)\Big(f_r(d_H(z)-1)-f_r(d_H(z))\Big)\\
&=&
\sum_{\{ y,z\}\in {\cal E}}\Big(h_1(y)h_2(z)+h_1(z)h_2(y)\Big)\mbox{, where}\\
h_1(x)&=&(r-1)d_H(x)-1\mbox{ and }\\
h_2(x)&=&f_r(d_H(x)-1)-f_r(d_H(x)).
\end{eqnarray*}
If $d_H(y)\geq d_H(z)$, then $h_1(y)\geq h_1(z)$ and, by Lemma \ref{lemma2}, $h_2(z)\geq h_2(y)$,
which implies
$$\Big(h_1(y)-h_1(z)\Big)\Big(h_2(z)-h_2(y))\Big)\geq 0.$$
Therefore,
$h_1(y)h_2(z)+h_1(z)h_2(y)\geq h_1(y)h_2(y)+h_1(z)h_2(z)$.

Since, for every vertex $y$ of $H$, there are exactly $(r-1)d_H(y)$ many vertices $z$ of $H$ with $\{ y,z\}\in {\cal E}$,
we obtain
\begin{eqnarray*}
P_2&\geq &
\sum_{\{ y,z\}\in {\cal E}}\Big(h_1(y)h_2(z)+h_1(z)h_2(y)\Big)\\
&\geq &
\sum_{\{ y,z\}\in {\cal E}}\Big(h_1(y)h_2(y)+h_1(z)h_2(z)\Big)\\
&=&
\sum_{x\in V(H)}(r-1)d_H(x)h_1(x)h_2(x)\\
&=&\sum_{x\in V(H)}(r-1)d_H(x)\Big((r-1)d_H(x)-1\Big)\Big(f_r(d_H(x)-1)-f_r(d_H(x))\Big).
\end{eqnarray*}
Combining these estimates, we see that 
\begin{eqnarray*}
P &=& P_1+P_2\\
&\geq& \sum_{x\in V(H)}
\Big((r-1)-(r-1)(d_H(x)+1)f_r(d_H(x))\\
&&\hspace{2cm}+(r-1)d_H(x)\Big((r-1)d_H(x)-1\Big)\Big(f_r(d_H(x)-1)-f_r(d_H(x))\Big)\Big),
\end{eqnarray*}
which is $0$ by the definition of $f_r$.
This completes the proof.
$\Box$

\bigskip

\noindent It seems a challenging task to extend the presented results to non-uniform and/or non-linear triangle-free hypergraphs.

\end{document}